\documentclass[12pt]{article}
\usepackage{amssymb,amsmath,amsfonts,latexsym,pap_n}

\newtheorem{Theorem}{Theorem}
\newtheorem{Definition}{Definition}

\newtheorem{Lemma}{Lemma}

\newtheorem{Corollary}{Corollary}

\newtheorem{Remark}{Remark}

\def\sqr#1#2{{\vcenter{\vbox{\hrule  height.#2pt
        \hbox{\vrule width.#2pt height#1pt \kern#1pt \vrule width.#2pt}
        \hrule height.#2pt}}}}
\def\Box{\sqr66}   
\let\epsilon=\varepsilon
\def\){ \right) }
\def\({ \left( }
\def\[{ \left[ }
\def\]{ \right] }
\def\<{ \langle }
\def\>{ \rangle }
\let\ljunk=\{
\let\rjunk=\}
\def\{{\left\ljunk}
\def\}{\right\rjunk}
\def\p{\partial}

\newcommand{\R}{{\mathbf R}}

\begin{document} 

\title{ Algebraic Aspects of the Theory of Product Structures in
Complex Cobordism} 
\author{B. Botvinnik, V. Buchstaber,
S. Novikov, S. Yuzvinsky}
\maketitle
\vspace{-1.5em}
\begin{abstract}
We address the general classification problem of all stable associative
product structures in the complex cobordism theory. We show how to
reduce this problem to the algebraic one in terms of the Hopf algebra
$S$ (the Landweber-Novikov algebra) acting on its dual Hopf algebra
$S^*$ with a distinguished ``topologically integral'' part $\Lambda$
that coincids with the coefficient ring of the complex cobordism.
We describe the formal group and its logarithm in terms of
representations of $S$. We introduce one-dimensional
representations of a Hopf algebra. We give series of examples of such
representations motivated by well-known topological and algebraic
results. We define and study the divided difference operators on an
integral domain.
 We discuss certain
important examples of such operators arising from analysis,
representation theory, and noncommutative algebra. We give a special
attention to the division operators by a noninvertible element of
a ring. We give new general constructions of associative product
structures (not necessarily commutative) using the divided
difference operators. As application, we describe new classes of
associative products in the complex cobordism theory.
\end{abstract}
\vspace{-1.5em}
\tableofcontents
\vspace{-.5em}
\section{Introduction. Complex cobordism theory and algebra}
Complex cobordism (bordism) functor $U^*(X)$, $U^*(X,Y)$ is a well-known
generalized cohomology (homology) theory with dual
elements $u\in U_k(X)$ being represented by maps $f: M^k \rightarrow
X$ of closed $U$-manifolds to $X$ (a complex structure is given
in the normal bundle of some embedding $M^k\subset \R^{n+k}$, see
\cite{1}). Algebraically, $U^*(X,Y)$ is a commutative and associative
$\Bbb Z$-graded ring with a unit $1\in U^*(pt)=\Omega^*_U = \Lambda$ (for
any pairs $(X,Y)$). According to the results of Milnor and Novikov (see
\cite{2}-\cite{4}), this ring is isomorphic to the polynomial ring
$\Omega^*_U=\Bbb Z[u_2,u_4,u_6,\ldots]$ with even-dimensional generators,
$\deg u_k = -2k$.
\vspace{2mm}

\noindent
{\bf Problem.} Give classification of all stable associative product
structures in the complex cobordism theory (on the category of finite
$CW$-complexes).
\vspace{2mm}

To make it precise, a ``product structure in the complex cobordism''
is a bilinear operation $u\circ v=\Phi(u\otimes v)$ where
$$
\Phi: U^*(X)\otimes U^*(X)\to U^*(X),
$$
defined for all spaces $X$, all $u,v\in U^*(X)$, and
commuting with continuous maps 
$f: X\to Y$, i.e.
$$
f^*(u)\circ f^*(v)=\Phi(f^*(u)\otimes f^*(v))=f^*\Phi(u\otimes v)
=f^*(u\circ v).
$$
A product structure is {\it stable} if it commutes with the
suspension isomorphism $s$.  Recall that there is a distinguished
element $\mu\in U^1(S^1)$ such that the suspension isomorphism $s:
U^*(X)\to\widetilde U^*(\Sigma X)$ is given by $u\to\mu u=s(u)$.  Here
$S^1$ is a circle, $\Sigma X=S^1\wedge X$ is a suspension over
$X$, and $\widetilde U^*(\Sigma X)=U^*(\Sigma X,pt)$.  Here and below
the symbol $\wedge$ denotes the smash-product in the category of
spaces with a base point, i.e. $X\wedge Y=X\times Y/X\vee Y$.  Thus a
new product structure $\Phi(u\otimes v)=u\circ v$ is stable if the
condition
\begin{equation}\label{eq1}
\mu(u\circ v)=(\mu u)\circ v=(-1)^{q_1}u\circ(\mu v)
\end{equation}
is satisfied. Here $q_1=\deg u$ and $uv$ denotes the usual standard
product in the complex cobordism. Clearly the standard product is stable.

We recall that the complex cobordism theory is a module over the
``Steenrod algebra'' of all operations $A^U$, i.e. the linear
operators $a: U^*(X)\to U^*(X)$ commuting with all continuous maps
and the suspension isomorphism. This algebra has been computed in \cite{5}
and consists of two parts:

1. The coefficient ring $\Lambda=\Omega^U\subset A^U$ that consists of the
   operations of multiplication by a ``scalar'' $\lambda\in\Lambda$:
   $u\to\lambda u$.

\let\varnothing=\emptyset 

2. The Landweber-Novikov algebra $S$ (see \cite{5}, \cite{6}). This is
   a Hopf algebra over $\Bbb Z$, it is $\Bbb Z^+$-graded,
   and $S=\sum_{j\ge0}S^j$ with a basis $s_w\in S$ where $S^0=\Bbb
   Z$. Here $w=(k_1,\dots,k_l)$ is either an unordered collection of
   positive integers $k_j\in\Bbb Z^+$, or the empty set with
   $s_\varnothing=s_0=1\in S^0$ and $\deg w=2\sum k_q$.

The coproduct (or diagonal) in  $S$ is given by
\begin{equation}\label{eq2}
\Delta s_w=\sum_{(w',w'')=w}s_{w'}\otimes s_{w''}.
\end{equation}
Furthermore, the Hopf algebra $S$ acts on all $U^*(X)$ (in particular
on the coefficient ring $\Lambda=U^*(pt)= \Omega^U$) providing
$U^*(X)$ with a structure of ``Milnor module'':
\begin{equation}\label{eq3}
s_w(uv)=\sum_{(w',w'')=w}s_{w'}(u)s_{w''}(v).
\end{equation}
The elements $s_{(n)}$ are {\it primitive} in $S$ and they generate a
Lie algebra; indeed, one has the identities
$[s_{(n)},s_{(m)}]=(m-n)s_{(m+n)}$.

An element $x\in U^2(X)$ is a {\it geometric cobordism element} if the
following conditions are satisfied
\begin{equation}\label{eq4}
s_{(k)}(x)=x^{k+1} 
\ \ \ \mbox{and} \ \ \ 
s_w(x)=0 
\ \ \ \mbox{if} \ \ \
w=(k_1,\dots,k_l), 
\ \ \ 
l>1.
\end{equation}
There is a natural representation 
$$
s\to s(x_1x_2\cdots x_N), 
\ \ \
s\in S,
$$
of the algebra $S$ (and the whole algebra $A^U$) on the geometric cobordism
elements. This representation is {\it asymptotically faithful} in the
sense that for any element $s\in S$ there exists an integer $N>m$ so
that $s(x_1x_2\cdots x_N)=0$ is equivalent to $s=0$. Recall also that
the representation $\lambda\to s\lambda$ of $S$ on the coefficient
ring $\Lambda$ is faithful. The latter was computed in \cite{5}.  The
paper \cite{7} gives interesting interpretation of this representation
in terms of differential operators on infinite dimensional Lie group
of formal derivations on the real line. This example leads to the
theory of operators on quantum groups (see \cite{8}, \cite{9}).
The algebra $A^U$ of operations in complex cobordism may be described
as the completion $A^U=(\Lambda S)^\wedge$, i.e. its elements are the series
$$
\sum^\infty_{j=0}\lambda_js_{w_j}\in A^U, 
\ \ \ 
\lambda_j\in\Lambda, 
\ \ \ 
s_{w_j}\in S,
$$
where $\deg w_j\to\infty$ while $j\to\infty$. The commutation relations in $A^U$
are determined by the representation of the algebra $S$ on the
coefficient ring $\Lambda$. 

Let $S^*$ be the Hopf algebra dual to $S$. Then there is the embedding
$\Lambda\subset S^*$ given by $\lambda(s_{w})= \varepsilon
s_w(\lambda)$ where $\lambda\in \Lambda$ and $\varepsilon
:\Lambda\to\Bbb Z$ is the augmentation ($\varepsilon(1)=1$ and
$\varepsilon(\lambda)=0$ if $\deg\lambda<0$).  Clearly $\Lambda\otimes
Q=S^*\otimes Q$ and $S^*$ is the polynomial ring over $\Bbb Z$ on
generators $s^*_k$ such that $(s^*_k,s_w)=1$ if $w=(k)$ and
$(s^*_k,s_w)=0$ otherwise. Below we give an algebraic description of
the subring $\Lambda\subset S^*$.
\begin{Lemma}\label{lemma1}
Any stable product structure $\Phi$ in the complex cobordism is given
by a formal series
\begin{equation}\label{eq5}
\widehat\Phi=\sum\lambda_{ij}s_{w_i}\otimes s_{w_j}
\end{equation}
where $\lambda_{ij}\in\Lambda$, $\deg w_i\to\infty$, $\deg
w_j\to\infty$, while $i,j\to\infty$, so that the product structure
$\Phi$ is given by
\begin{equation}\label{eq6}
u\circ v=\Phi(u\otimes v)=\sum\lambda_{ij}s_{w_i}(u)s_{w_j}(v)
\end{equation}
for any $u,v\in U^*(X)$. 

Furthermore, the coefficients of the series $\widehat\Phi$ are
uniquely determined by the products $u\circ v$ for
$u,v\in\Lambda=U^*(pt)=\Omega^U$.
\end{Lemma}
\begin{Proof}
Let $x(N)\in U^{2N}(MU_N)$ be a canonical Thom class corresponding to
the universal complex bundle $\eta_N \to BU(N)$ where $\dim
\eta_N=N$ and $MU_N$ is the Thom space of this bundle. 

It is well-known that $U^*(MU_N,pt)$ is one-dimensional $A^U$-module
generated by $x(N)$ and the elements $s_wx(N)$ with
$w=(k_1,\dots,k_l)$, $l\le N$, $x(N)=s_0x(N)$, form its basis over 
$\Lambda$. 

Let $\Phi$ be a stable product structure in the complex cobordism
theory. Then the products $x(N_1)\circ x(N_2)\in U^*(MU_{N_1}\times
MU_{N_2})$ are defined for any $N_1$, $N_2$. The bilinearity of $\Phi$
implies
$$
\Phi(x_{N_1},x_{N_2})=
x(N_1)\circ x(N_2)\in\widetilde U^*(MU_{N_1}\wedge MU_{N_2}).
$$
Thus this product may be uniquely written as a formal series
$$
x(N_1)\circ x(N_2)=\sum\lambda_{ij}s_{w_i}(x(N_1))s_{w_j}(x(N_2))
$$
where $w_i=(k_1',\dots,k_{l_1}')$, $l_1\le N_1$, and
$w_j=(k_1'',\dots,k_{l_2}'')$, $l_2\le N_2$. The coefficients
$\lambda_{ij}\in\Lambda$ do not depend on $N_1$ and $N_2$ since the
product is stable and commutes with continuous maps. We recall here
that the Thom spaces $MU_N$ for different $N$ are related to each
other by means of the canonical maps $e_{N_1,N_2}
:\Sigma^{2(N_2-N_1)}MU_{N_1}\to MU_{N_2}$ such that
$e^*_{N_1,N_2}x(N_2)=s^{2(N_2-N_1)}x(N_1)$.

Let $u\in U^{q_1}(X,Y)$ and $v\in U^{q_2}(X,Y)$ be any two
elements. Then there exist maps
$$
f_1 : \Sigma^{2N_1-q_1}X/Y\to MU_{N_1} 
\ \ \ \mbox{and} \ \ \ 
f_2 : \Sigma^{2N_2-q_2}X/Y\to MU_{N_2}
$$
representing $u$ and $v$ for some $N_1$ and $N_2$,
i.e. $f^*_1x(N_1)=s^{2N_1-q_1}u$, $f^*_2x(N_2)=s^{2N_2-q_2}v$. We
define a map
$$
f=f_1\wedge f_2 :\Sigma^{2(N_1+N_2)-(q_1+q_2)}X
\stackrel{f_1\wedge f_2}{\longrightarrow}
MU_{N_1}\wedge MU_{N_2}.
$$
Then we have
\begin{equation*}
\begin{split}
f^*(x(N_1)\circ x(N_2))
&=\sum\lambda_{ij}s_{w_i}(s^{2N_1-q_1}u)
s_{w_j}(s^{2N_2-q_2}v)
\\
&=(-1)^{q_1q_2}s^{2(N_1+N_2)-(q_1+q_2)}
\Bigl(\sum\lambda_{ij}s_{w_i}(u)s_{w_j}(v)\Bigr).
\end{split}
\end{equation*}
On the other hand,
\begin{equation*}
\begin{split}
f^*(x(N_1)\circ x(N_2))
&=(f^*_1x(N_1))\circ(f^*_2x(N_2))=(s^{2N_1-q_1}u)\circ(s^{2N_2-q_2}v)
\\
&=(-1)^{q_1q_2}s^{2(N_1+N_2)-(q_1+q_2)}u\circ v.
\end{split}
\end{equation*}
This proves the first statement of Lemma \ref{lemma1}. 

Now let $u\circ v$ be a product on the coefficient ring
$\Lambda=U^*(pt)$. We notice that if $u\in\Lambda^{-2n}$ and 
$v\in\Lambda^{-2m}$ then
\begin{equation}\label{eq7}
u\circ v=\Phi(u\otimes v)=
\sum_{\begin{array}{c}^{\deg w_i\le2n,} \\ ^{\deg w_j\le2m}\end{array}}
\lambda_{ij}s_{w_i}(u)s_{w_j}(v).
\end{equation}
In particular, for $\deg w=2n$ we have 
$$
s_w(u)=\varepsilon s_w(u)=\langle u,s_w\rangle
$$
where $\varepsilon :\Lambda\to Z$ is the augmentation.

Now we use 
(\ref{eq7}), 
the above embedding $\Lambda\subset S^*$, and the isomorphism
$\Lambda\otimes Q=S^*\otimes Q$ to recover the coefficients
$\lambda_{ij}$ by induction on $\deg w$. Since we use only the product
$u\circ v$ on $\Lambda$ and the action of algebra $S$ on
$\Lambda$, this completes the proof of the  lemma.  
\end{Proof}

Lemma \ref{lemma1} 
allows us to define a new product in the complex cobordism $U^*(X)$ as a
series $\widehat\Phi$ given by (\ref{eq5}). Below we study properties of
such products.

Now we describe this situation in terms of 
Hopf algebras.
Suppose $S$ is a Hopf algebra and $\Lambda$ is a Milnor module 
over $S$ (i.e. $S$ acts on $\Lambda$ and $s(uv)=\sum
s_i'(u)s_i''(v)$ where $\Delta(s)=\sum s_i'\otimes s_i''$). Then the new
product structure is defined by (\ref{eq6}).  This gives new setting
for the classification problem of product structures.

The most interesting case arising from the complex cobordism theory
is when the Milnor module is $S^*$, that is the Hopf algebra dual to
$S$ in the basis $s_w$ with $\Lambda\subset S^*$. The action of $S$
on $S^*$ is given by $s(u)=R^*_s(u)$ where $R_s$ is the 
{\sl right multiplication operator} (see \cite{7}, \cite{8}):
\begin{equation}\label{eq8}
R_s(s')=s's, \quad R_s: S\to S, \quad
R^*_s: S^*\to S^*, \quad (R^*_s(u),s')=(u,s's).
\end{equation}
Thus the classification problem of product
structures in complex cobordism $U^*(X)$ reduces to study of the single Hopf
algebra $S$ acting (by means of $R^*$) on its dual Hopf algebra $S^*$
with a distinguished ``topologically integral'' part
$\Lambda\subset S^*$.

How to describe the ring $\Lambda\subset S^*$ in algebraic
terms? For this we use geometric cobordism elements
$x\in U^2(X)$ defined above.  Consider the Milnor module
$S^*[[x_1,x_2]]$ of formal series on $x_1$, $x_2$ over $S^*$. 
Define an action of $S$ on this module algebraically, provided $x_1$
and $x_2$ are geometric cobordism elements, by
\begin{equation*}
\begin{array}{c}
s_{(k)}(x_j)=x^{k+1}_j, \qquad s_w(x_j)=0, \quad w=(k_1,\dots,k_l), \quad l>1,
\\
s(\lambda)=R^*_s(\lambda), \qquad (R^*_s(\lambda),s')=(\lambda,s's),
\\
s(uv)=\sum s_i'(u)s_i''(v), \qquad \Delta=\sum s_i'\otimes s_i''.
\end{array}
\end{equation*}
This determines the action of $S$ on $S^*[[x_1,x_2]]$.
\begin{Lemma}\label{lemma2}
There exists a unique series
\begin{equation}\label{eq9}
x=f(x_1,x_2)=x_1+x_2+\sum_{i,j\ge1}\alpha_{ij}x^i_1x^j_2, \qquad
\alpha_{ij}\in S^*,
\end{equation}
such that $x=f(x_1,x_2)$ is a geometric cobordism element provided
$x_1$ and $x_2$ are geometric cobordism elements. In particular, the
action of $S$ on $x$ is determined by {\rm (\ref{eq4})}.
\end{Lemma}
\begin{Proof}
We define the operator $S_t=\sum_{w\ge0}s_wt^w$ with
$w=(k_1,\dots,k_l)$, $t^w=t^{k_1}_1\dotsb t^{k_l}_l$, where $t_i$,
$i=1,2,\dots$, are algebraically independent elements. It is easy to
see that
$$
\Delta S_t=\sum_{w\ge0}(\Delta s_w)t^w=S_t\otimes S_t \quad
\mbox{or}
\quad
S_t(uv)=S_t(u)S_t(v).
$$
Let us compute $S_t(x)$ for any geometric cobordism element $x$. By
definition of geometric cobordism elements, $S_t(x)$ is equal to
certain series $\phi_t(x)=\sum^\infty_{k=0}x^{k+1}t_k$ with
$t_0=1$. Now we apply the operator $S_t$ to the both sides of
(\ref{eq9}). 
\begin{equation}\label{eq10}
S_t(x)=S_t(x_1)+S_t(x_2)+\sum_{i,j\ge1}S_t(\alpha_{ij})S_t(x_1)^iS_t(x_2)^j
\end{equation}
where $x$, $x_1$, $x_2$ are geometric cobordism elements, i.e.
$S_t(x)=\phi_t(x)$, $S_t(x_1)=\phi_t(x_1)$,
$S_t(x_2)=\phi_t(x_2)$. Now we apply the augmentation $\varepsilon$
to (\ref{eq10}): 
$$
\varepsilon : S^*\to\Bbb Z
$$
where $\varepsilon(s^*_n)=0$ if $n\ne0$, $S^*=\sum_{n\le0}S^*_n$, and
$\varepsilon(x_1)=x_1$, $\varepsilon(x_2)=x_2$. We obtain
$$
\begin{array}{l}
\varepsilon S_t(x)=\sum\varepsilon(x)^{k+1}t_k, \qquad
\varepsilon(x)=x_1+x_2,
\\
\varepsilon S_t(x_1)=\sum\varepsilon(x_1)^{k+1}t_k=\sum x^{k+1}_1t_k,
\\
\varepsilon S_t(x_2)=\sum\varepsilon(x_2)^{k+1}t_k=\sum x^{k+1}_2t_k.
\end{array}
$$
The following formula 
$$
\varepsilon S_t(\alpha_{ij})=\sum_w\varepsilon s_w(\alpha_{ij})t^w
=\sum_w(s_w,\alpha_{ij})t^w
$$
follows from definitions. It implies
$$
\varepsilon S_t(x)=\phi_t(x_1+x_2)=\phi_t(x_1)+\phi_t(x_2)+\sum_{i,j}\sum_w
(s_w,\alpha_{ij})t^w\phi_t(x_1)^i\phi_t(x_2)^j.
$$
Let $y_j=\phi_t(x_j)$. Then $x_j=\phi^{-1}_t(y_j)$ and 
$$
\varepsilon S_t(x)=\phi_t(\phi^{-1}_t(y_1)+\phi^{-1}_t(y_2))=\sum_{i,j}\sum_w
(s_w,\alpha_{ij})t^wy^i_1y^j_2.
$$
Compare the coefficients of the terms $y^i_1y^j_2$ in the left
and right hand sides. Since $(s_w,\alpha_{ij})$ are already known
on the left-hand side, $\alpha_{ij}$ are uniquely
determined as elements of $S^*$, i.e. as linear forms on $S$. This
proves the uniqueness.

The existence of the series (\ref{eq9}) follows from well-known
results on the complex cobordism theory (see \cite{5}). An algebraic proof
of its existence follows from Lemma \ref{lemma3} below. This completes
the proof of Lemma \ref{lemma2}.
\end{Proof}
Lemma \ref{lemma2} implies the following.
\begin{Corollary}\label{corollary1}
The series $x=f(x_1,x_2)$ determines a formal group on the set of
geometric cobordism elements (see {\rm \cite{9}}).
\end{Corollary}
\begin{Definition}\label{definition1}{\rm 
The subring $\Lambda \subset S^*$ generated over $\Bbb Z$ by the
elements $\alpha_{ij}\in S^*$ is called the complex cobordism ring.
(It is well-known that $\Lambda$ is a polynomial ring, see above.)}
\end{Definition}
\begin{Remark}{\rm 
A formal group in complex cobordism theory
was defined in \cite{5} geometrically. It
coincides with the one given by the series
(\ref{eq9}). 
We emphasize that we give here entirely algebraic definition of this
formal group and the complex cobordism ring $\Lambda$ using only
the Hopf algebra $S$.}
\end{Remark}
The formal group $f(x_1,x_2)$ over the ring $\Lambda\otimes\Bbb Q$ is
given by
$$
f(x_1,x_2)=g^{-1}(g(x_1)+g(x_2)),
$$
where $g^{-1}(g(x))=x$ and $g(x)=x+\sum b_ix^{i+1}$,
$b_i\in\Lambda\otimes\Bbb Q$. The series $g(x)=g_f(x)$ 
is called the {\it logarithm of the formal group $f(x_1,x_2)$}.
We would like to describe the logarithm $g(x)$ again by means of 
the Hopf algebra $S$.
\begin{Lemma}\label{lemma3}
There exists a unique series
$$
g(x)=x+\sum^\infty_{i=1}b_ix^{i+1}, \qquad b_i\in S^*,
$$
such that $s_wg(x)=0$ 
for any geometric cobordism element $x$ and for all $w$ with $\deg w>0$. 
This series is a logarithm of the formal
group and satisfies the condition
$$
x=g(x)+\sum^\infty_{k=1}s^*_kg(x)^{k+1},
$$
i.e. $g^{-1}(t)=t+\sum^\infty_{k=1}s^*_kt^{k+1}$ where $s^*_k$ are
multiplicative generators of $S^*$ dual to $s_{(k)}\in S$ (with
respect to the $\Bbb Z$-basis $\{s_w\}$).
\end{Lemma}
\begin{Proof}
Let $x$ be a geometric cobordism element and let $g(x)=x+\sum
b_ix^{i+1}$ ($b_i\in S^*$) be a series such that $s_wg(x)=0$ for
all $w$ with $w\ne\emptyset$. Then $S_tg(x)=g(x)$ where $S_t=\sum
s_wt^w$. It implies that $x=g^{-1}(g(x))$. Thus
$\phi_t(x)=S_tx=S_t[g^{-1}](S_tg(x))=S_t[g^{-1}](g(x))$ where if
$g^{-1}(t)=t+\sum a_it^{i+1}$, $a_i\in S^*$ then $S_t[g^{-1}]=t+\sum
S_t(a_i)t^{i+1}$. Now let $x_1$ and $x_2$ be geometric cobordism
elements. Consider the series
$$
F(x_1,x_2)=g^{-1}(g(x_1)+g(x_2))=x_1+x_2+\sum\alpha_{ij}'x^i_1x^j_2.
$$
We have
$$
g(x_1)+g(x_2)=g(F(x_1,x_2)).
$$
Then 
$$
\begin{array}{rl}
S_tF(x_1,x_2)&=S_t[g^{-1}](S_tg(x_1)+S_tg(x_2))
\\
&=S_t[g^{-1}](g(x_1)+g(x_2))=S_t[g^{-1}](g(F(x_1,x_2))).
\end{array}
$$
We use the identity $S_t[g^{-1}](g(x))=\phi_t(x)$ to obtain
$$
S_tF(x_1,x_2)=\phi_t(F(x_1,x_2)).
$$
We see that the series $F(x_1,x_2)$ is a geometric cobordism element
whence according to Lemma \ref{lemma2}, it coincides with the series
(\ref{eq9}).  
Thus if a series $g(x)$ satisfies the conditions of Lemma \ref{lemma3}
then it is unique and coincides with the logarithm of the formal group
(\ref{eq9}).

Let us give an algebraic proof of the existence of such series. First, we
consider $S^*[[t]]$ as a Milnor module with the following action of
the algebra $S$:
$$
s_0t=t, \quad s_wt=0, \quad \deg w>0 \quad 
\mbox{and}
\quad
s_w(\lambda)=R^*_s(\lambda), \quad \lambda\in S^*.
$$
We show that the series 
$$
x=\gamma(t)=t+\sum s^*_{(n)}t^{n+1}\in S^*[[t]]
$$
is a geometric cobordism element. It is enough to prove that
$$
s_w(s^*_{(n)})= 
\{\begin{array}{ll}
0,& 
\mbox{if} \ w\ne(k),
\\
\gamma_{n,k},&
\mbox{if} \ w=(k),
\end{array}\right.
$$
where $\gamma_{n,k}$ are the coefficients of $t^{n+1}$ in the series
$\gamma(t)^{k+1}$. To complete the argument we use the formula
$$
\langle s_{w_2}(s^*_n),s_{w_1}\rangle=\langle s^*_n,s_{w_1}s_{w_2}\rangle
$$
and the fact that the representation of $S$ on a product of geometric
cobordism elements is asymptotically faithful. 

Now let $g(x)\in S^*[[t]]$ be such that $g(x)=g(\gamma(t))=t$. Then
$x$ is a geometric cobordism element and by construction $s_w(t)=0$
if $\deg w>0$. Thus the series $g(x)$ satisfies all conditions given
in Lemma \ref{lemma3}. This completes the proof.
\end{Proof}
We remark that the result of Lemma \ref{lemma3} has been proven first by
means of the complex cobordism theory (see \cite{10}). 
\begin{Definition}
{\rm We call a representation of a Hopf algebra $S$ on a Milnor module $P$
one-dimensional if $P$ is the polynomial algebra (or the algebra of formal
series) with one generator $u\in P$ over some coefficient ring. 
The algebra $S$ may act nontrivially on the coefficient ring
(which is $\Bbb Z$, $\Bbb Q$, $\Lambda$ in our examples, or some other
Milnor submodule in $S^*\otimes\Bbb Q$). A {\it height of such
representation} is a minimal number $k$ so that  $u^{k+1}=0$.}
\end{Definition}
For instance, $P$ may be a polynomial algebra with one generator $u\in P$ 
over the ring 
$\Lambda\subset S^*$ where $S$ acts as it was described above.
\vspace{2mm}

\noindent
{\bf Examples.} If $k=1$ we can let $P=U^*(S^{2q})$ for
any $q$.  Already for $k=2$, one-dimensional representations over $\Bbb Z$
are known only for 
$\dim u=2,4$ or $8$. For $k\geq 3$, we know examples from topology only for the
cases
$\dim u=2$ or $4$ and $P=U^*(\Bbb CP^k)$ or
$P=U^*(\Bbb HP^k)$ respectively. The case $k=\infty$ is of special
interest. For instance, over $\Bbb Q$ there are one-dimensional
representations of $S$ on $P=U^*(\Omega S^{2q+1})\otimes\Bbb
Q=\Omega_U\otimes\Bbb Q[[u]]$ with $\dim u_q=2q$ for all $q$.

Consider the case $q=1$ in more detail. A generator of the cohomology
group $H^2(\Omega S^3,\Bbb Z)=\Bbb Z$ is represented by a map $\phi
:\Omega S^3\to\Bbb CP^\infty$. The map $\phi$ gives a geometric
cobordism element $u_1=\phi^*(u)\in U^2(\Omega S^3)$ where $u\in
U^2(\Bbb CP^\infty)$ is a canonical geometric cobordism element. We
have the isomorphism $U^*(\Omega S^3)\otimes\Bbb Q=\Omega_U\otimes\Bbb
Q[[u_1]]$. It is well-known that there exists a map $\psi :
\Sigma\Omega S^3\to S^3$ giving an element $\beta\in U^2(\Omega
S^3)$ such that $s\beta=\psi^*\alpha$ where $\alpha$ is a generator
of the group $U^3(S^3)=\Bbb Z$ and $s$ is the suspension
isomorphism. All operations $s_w$ are stable whence $s_w\alpha=0$ for $\deg w>0$
implies
$s_w\beta=0$ for $\deg w>0$.
The element $\beta\in U^*(\Omega S^3)\otimes\Bbb Q$ is given by a
series $\beta=g(u_1)$ where $u_1$ is a geometric cobordism
element and $g(u_1)=u_1+\dotsb$. Thus Lemma \ref{lemma3} implies that 
$g(u_1)$ coincides with the logarithm of the above formal group and
$g^{-1}(\beta)=\beta+\sum_{k\ge1}s^*_k\beta^{k+1}=u_1\in U^2(\Omega
S^3)$. We emphasize that the elements $s^*_k\in S^*$ do not belong to
the subring $\Lambda\subset S^*$, only their multiples are:
$(k+1)!s^*_k\in\Lambda$ (see \cite{10}).

To describe one-dimensional representations of $S$ on $P=U^*(\Omega
S^{2q+1})\otimes\Bbb Q$ (in the case $q>1$) one should use the canonical
maps
$$
\phi_q : \Omega S^{2q+1}\to\Omega SU(2q+1)\to BU
$$
representing the corresponding generators of the homotopy group
$\pi_{2q+1}(SU(2q+1))=\Bbb Z$. Here $SU(2q+1)$ is a special unitary
group and $BU$ is the infinite dimensional complex Grassmannian. In order to
construct such maps one uses that the Bott periodicity implies
$BU\approx\Omega SU(\infty)$.

There is an important example of one-dimensional (graded) $\Bbb
Z$-representation of the Hopf algebra $S$ given by a geometric
cobordism elements $u=x$ (as above) with $\dim u=2$. One more example
of one-dimensional $\Lambda$-representation of $S$ is given by the
elements $u=x\bar x$, $\dim u =4$, where $\bar x$ is the element
inverse to $x$ in the formal group. Here $s_w(u)\in \Lambda[[u]]$.
Such elements $u$ generate a two-valued formal group (this was first
observed in \cite{11}).  The algebraic theory of two-valued formal
groups was developed in \cite{12}. The case $\dim u=8$, $k=2$ is also
very interesting. Here $P=U^*(CaP^2)$ where $CaP^2$ is the projective
plane over the Kelley numbers.  We suspect that there are no
one-dimensional representations of $S$ for $k\geq 3$ with $\dim
u>4$. It would be useful to prove this.
\vspace{2mm}

\noindent
{\bf Problem.} Give classification of one-dimensional $\Bbb
Z$-representations and $\Lambda$-representations of the Hopf algebra
$S$. The representations should preserve grading under the
automorphisms
$$
u=\Psi(v)=v+\sum_{i\ge1}\lambda_iv^{i+1}
$$
where $\lambda_i\in\Bbb Z$ (or $\Bbb Q$), or $\lambda_i\in\Lambda$ (or
$\lambda_i\in\Lambda\otimes\Bbb Q=S^*\otimes\Bbb Q$). 
\vspace{2mm}

Some examples of such representations over certain subrings of
$\Lambda\otimes\Bbb Q$ may be found in \cite{13}. They are given by
the elements
$u_n=x[x]_{\varepsilon_n}\dotsb[x]_{\varepsilon^{n-1}_n}$ with $\deg
w =2n$. Here $\varepsilon_n$ is the $n$th root of unity and
$[x]_{\varepsilon_n}=g^{-1}(\varepsilon_ng(x))$ is a geometric
cobordism element that is nothing but the $\varepsilon_n$-th power
in the formal group $f(x_1,x_2)$ of the canonical geometric cobordism
element $x$. According to \cite{13}, $u_n\in\Lambda_{\{n\}}[[x]]$
where $\Lambda_{\{n\}}=\Lambda\otimes\Bbb Z_{\{n\}}$. Here
$$
\Bbb Z_{\{n\}}= \{ \frac{m}{d} \in \Bbb Q \ \left| 
\begin{array}{c} (d,p) =1 \ \mbox{if the equation $x^n=1$ has precisely $n$}
\\
\mbox{distinct solutions in the ring of $p$-adic
numbers}
\end{array}\right.
\}
\subset\Bbb Q .
$$ 
For instance, one has $\Bbb Z_{\{2\}}=\Bbb Z$. This
case gives $\Lambda$-representations for $u=x\bar x$ (see above).
These $\Lambda$-representations are related to the theory of
two-valued formal groups. Furthermore, one-dimensional
$\Lambda_{\{n\}}$-representations of $S$ on
$P=\Lambda_{\{n\}}[[u_n]]$ with $\deg u_n=2n$ are related to
$n$-valued formal groups. The existence of $n$-valued formal groups
was observed first in \cite{11}. First results on the algebraic theory
of these groups are given in \cite{14}. The contemporary state of affairs in
the theory of multi-valued formal groups may be found in \cite{15}. 

There is another useful example of one-dimensional $\Bbb
Z$-representation of the Hopf algebra $S$ with $\dim u=2$. This
representation is not equivalent to any action on a geometric cobordism
element (see \cite{7}). Here is the action:
\begin{equation}\label{eq11}
s_{(1)}(u)=u^2, \qquad s_{(n)}(u)=0, \quad n\ge2 .
\end{equation}
We recall that one has $s_{(n)}(u)=u^{n+1}$ on a geometric cobordism
element $u$ (see (\ref{eq4}) above).  It is easy to see that the
representation (\ref{eq11}) of $S$ on $P=\Bbb Z[u]$ is
given by the differential operators
$$
s_{(1)}\to u^2\frac d{du}\,, \qquad s_{(n)}\to0, \quad n\ge2.
$$
Furthermore, the image of the Hopf algebra $S$ (in the algebra of
differential operators) is generated by $\frac{1}{n!}(u^2d/du)^n$.

There is an interesting class of such representations where $P$ is a
polynomial algebra on one generator with negative degree, i.e. $P=\Bbb Z[u]$
($\Bbb Q[u]$) or $P=\Lambda[u]$ ($\Lambda\otimes\Bbb Q[u]$). Here
$\dim u=-2q<0$.
\vspace{2mm}

\noindent
{\bf Example.} Let $u\in \Lambda=\Omega^U$ be an element with
$\dim u= -2$ or $-4$. If $u= \in\Lambda^{-2}$ let $u=[\Bbb CP^1]$
and if $u\in\Lambda^{-4}$ let $u=3[\Bbb CP^1]^2-4[\Bbb CP^2]$. Notice
that in the latter case $s_{(1)}(u)=0$. In both cases the action of
$S$ on $u$ is determined by the operations $s_w(u)\in\Bbb Z$ with
$\deg w=2$. Thus we have one-dimensional $\Bbb
Z$-representations of the algebra $S$.

\section{Division operators and product operators}
Here we develop an algebraic machinery which allows us to produce
examples of stable product structures in the complex cobordism theory.
We will carry on all computations  in the category of modules over
some commutative associative ring $K$ (the ``scalars'') with the
unit $1\in K$. Furthermore, all modules $R$ in our category are also
commutative and associative rings with the unit $1\in K\subset R$.
We also assume that $K$ and $R$ are integral domains. All operators
below are assumed to be $K$-linear.
\begin{Definition}\label{def3}{\rm 
A linear operator $\partial : R\to R$ is called a {\it divided
difference operator} if it is not identically trivial and satisfies
the following identity
\begin{equation}\label{eq12}
\partial(xy)=(\partial x)y+x(\partial y)-\alpha(\partial x)
(\partial y)
\end{equation}
for all $x,y\in R$, where $\alpha\in R$ 
is not invertible in $R$.}
\end{Definition}
It follows from (\ref{eq12}) that $\partial(1)=0$ for any divided
difference operator $\partial$ (under our assumptions).

An operator $\pi : R\to R$ is {\it multiplicative} if
$\pi(xy)=\pi(x)\pi(y)$.  
\begin{Lemma}\label{lemma4} 
An operator $\partial$ is a divided difference operator if and only if
the operator $\pi=1-\alpha\partial$ is multiplicative.
\end{Lemma}
\begin{Proof}
Let $\partial$ be a divided difference operator. Then the definition
gives:
$$
\begin{array}{rl}
\pi(xy)&= \ xy-\alpha\partial(xy)
=xy-\alpha(\partial x)y-\alpha x(\partial y)+\alpha^2(\partial x)(\partial y)
\\
&= \ (x-\alpha\partial x)(y-\alpha\partial y)=\pi(x)\pi(y).
\end{array}
$$
It is easy to see that the converse statement holds as well. This
proves Lemma \ref{lemma4}.
\end{Proof}
Lemma \ref{lemma4} motivates the term ``divided difference operator''.
\begin{Lemma}\label{lemma5} 
Let $\partial$ be a divided difference operator such that
$\partial^2=\gamma\partial$ with $\gamma\in\Bbb R$. Then 
$(1-\alpha\gamma)\partial(\alpha)=2-\alpha\gamma$. In particular,
$\partial(\alpha)=2$ if $\partial^2=0$.
\end{Lemma}
\begin{Proof}
Indeed, the definition gives:
$$
\begin{array}{rl}
\partial^2(xy)
&=\partial\bigl((\partial x)y+x(\partial y)-
\alpha(\partial x)(\partial y)\bigr)
\\
&=(\partial^2x)y+(\partial x)(\partial y)-\alpha(\partial^2x)(\partial y)
+(\partial x)(\partial y)+x(\partial^2y)-\alpha(\partial x)(\partial^2y)
\\
&\qquad-(\partial\alpha)(\partial x)(\partial y)
- \alpha\partial\bigl((\partial x)(\partial y)\bigr)
+\alpha(\partial\alpha)\partial\bigl((\partial x)(\partial y)\bigr).
\end{array}
$$
The condition $\partial^2=\gamma\partial$ implies the identity
$$
(1-\alpha\gamma)[(1-\alpha\gamma)\partial(\alpha)-(2-\alpha\gamma)]
\partial x\partial y=0.
$$
We recall that $\partial$ is not identically trivial, $R$ is an
integral domain, and $\alpha$ is not invertible in $R$. This proves
the result.
\end{Proof}
\begin{Lemma}\label{lemma6}
A divided difference operator $\partial$ satisfies the condition
$\partial^2=\gamma\partial$ if and only if $\pi^2=1$ and 
$\pi(\alpha)=-\alpha/(1-\alpha\gamma)$ where $\pi=1-\alpha\partial$.
\end{Lemma}
\begin{Proof}
Lemma \ref{lemma4} asserts that  $\pi$ is a multiplicative operator. Thus
\begin{equation}\label{eq13}
\begin{array}{rl}
\pi^2x&=\pi(x-\alpha\partial x)=\pi(x)-\pi(\alpha)\pi(\partial x)
\\
&=x-\alpha\partial x-\pi(\alpha)\partial x+\alpha\pi(\alpha)\partial^2x
=x-(\alpha+\pi(\alpha))\partial x+\alpha\pi(\alpha)\partial^2x.
\end{array}
\end{equation}
Let  $\partial^2=\gamma\partial$. Then Lemma \ref{lemma5} gives 
$\partial(\alpha)=(2-\alpha\gamma)/(1-\alpha\gamma)$ whence 
$\pi(\alpha)=\alpha-\alpha\partial\alpha=-\alpha/(1-\alpha\gamma)$.
Thus the condition $\partial^2=\gamma\partial$ implies that
$\pi^2x=x$ for all $x$. Conversely, the conditions  $\pi^2=1$ and
$\pi(\alpha)=-\alpha/(1-\alpha\gamma)$ together with (\ref{eq13}) imply
$\partial^2=\gamma\partial$. This proves Lemma \ref{lemma6}.
\end{Proof}
The constructions given here suggest an interpretation of $R$ as a
space of functions on some space and operator $\pi$ as
a translation operator of that space or another homeomorphism of this space.
Of course
this interpretation makes sense if the operator $\pi$ is invertible.
In important examples, $\pi$ turns out to be a projector, i.e. $\pi^2=\pi$.
In order to construct multiplicative projectors,
the ``division operators'' (i.e. operators dividing by a
given scalar) were used in \cite[p.887]{5}. These operators are
defined by the identity (\ref{eq12}) and the requirement
$\partial(\alpha)=1$.  The paper \cite{16} gives a classification of
such operators, in particular, it was shown (see Lemma \ref{lemma12}
below) that the well-known Adams-Quillen projectors are compositions
of the Novikov's division operators by some scalars $\alpha$ (under
a special choice of these scalars).
\begin{Lemma}\label{lemma7}
A divided difference operator $\p$ is a division operator by a scalar
$\alpha$ if and only if the composition $\partial\alpha=\partial\circ\alpha$ 
is the identity
operator, in particular,
$$
\pi(x)=(\partial\alpha-\alpha\partial)(x)=[\partial,\alpha](x).
$$
\end{Lemma}
\begin{Proof}
The definitions imply
$$
\partial\alpha(x)=\partial(\alpha x)=(\partial(\alpha))x
+\alpha(\partial x)-\alpha(\partial(\alpha))(\partial x).
$$
If $\partial(\alpha)=1$ then $\partial\alpha x=x+\alpha(\partial
x)-\alpha(\partial x)=x$. Conversely, let $\partial(\alpha)=1$, i.e.
$\partial(\alpha x)=x$. This implies
$$
\partial(\alpha x)=(\partial(\alpha))x+\alpha(\partial x)-\alpha
\partial(\alpha)(\partial x)=x.
$$
Thus
$$
(1-\partial(\alpha))x=\alpha(\partial x)(1-\partial(\alpha)),
$$
or 
\begin{equation}\label{eq14}
(1-\partial(\alpha))[x-\alpha(\partial x)]=0.
\end{equation}
Now we let $x=1$ in (\ref{eq14}). Then  $\partial(1)=0$
gives $\partial(\alpha)=1$. The result follows.
\end{Proof}
We assume now that the ring $R$ satisfies the following condition. If
$\alpha\in R$ is a noninvertible element then $\bigcap_n(\alpha^n
R)=0$ where $(\alpha R)\subset R$ is the principal ideal generated by $\alpha$.
This assumption holds in all our examples.
\begin{Lemma}\label{lemma8}
Let $\partial$ be a divided difference operator. Then the kernel of
$\pi=1-\alpha\partial$ is nontrivial 
if and only if $\partial$ is a division operator.
\end{Lemma}
\begin{Proof}
Let $\partial$ be a division operator, i.e. $\partial(\alpha)=1$. Then 
$\pi(\alpha)=\alpha-\alpha\partial(\alpha)=0$. Conversely, let 
 $\partial(\alpha)\ne1$ and there exists an element $x\ne0$ such that
$\pi(x)=0$. Then we have
$$
x=\alpha\partial x.
$$
Appling $\partial$ we have
$$
\partial x=\partial(\alpha)\partial x+\alpha\partial^2x
-\alpha\partial(\alpha)\partial^2x.
$$
Thus 
$$
(1-\partial(\alpha))\partial x=\alpha(1-\partial(\alpha))\partial^2x.
$$
Now the condition $\partial(\alpha)\ne1$ implies
$$
\partial x=\alpha\partial^2x.
$$
We use induction to obtain the identity 
$$
\partial^kx=\alpha\partial^{k+1}x
$$
for any $k\ge0$. This implies 
$$
x=\alpha^k\partial^kx
$$
for all $k\ge0$, i.e. $x\in\bigcap_k(\alpha^kR)$. However, it is
possible only if $x=0$ because of the assumptions on $R$.  This
proves Lemma \ref{lemma8}.
\end{Proof}
\begin{Corollary}\label{cor2}
Let $\partial$ be a divided difference operator. Then the multiplicative
operator $\pi=1-\alpha\partial$ is a projector ($\pi^2=\pi$)
if and only if $\partial$ is a division operator by the element $\alpha$.
\end{Corollary}
\begin{Proof}
Let $\partial$ be a division operator by the element $\alpha$.
Then $\pi(\alpha)=0$ and $\pi^2x=\pi(x-\alpha\partial x)=\pi x$.
Conversely, let $\pi^2=\pi$. Then $\pi(\pi x-x)=0$ for any $x$ and
Lemma \ref{lemma8} gives that  $\pi x=x$. This contradicts 
$\alpha\ne0$. This proves the result.
\end{Proof}
\begin{Lemma}\label{lemma9}
 Let $\alpha\in R$ and let $\Pi_\alpha(R)$ be the
 set of multiplicative operators $\pi\in\Pi_\alpha(R)$ such that
 $\pi=1-\alpha\partial$. Then the set $\Pi_\alpha(R)$ is a
 subsemigroup in the semigroup $\Pi(R)$ of all multiplicative
 operators $R\to R$ where operation in $\Pi(R)$ is given by
 $(\pi_1\circ\pi_2)(u)=\pi_1(\pi_2(u))$. 

Let $\pi_1, \pi_2\in \Pi_\alpha(R)$ and $\partial_1$, $\partial_2$ be
the corresponding divided difference operators. Then $\partial_1$ and
$\partial_2$ satisfy
\begin{equation}\label{eq15}
\partial_1\circ\partial_2=\partial_1+\partial_2-\partial_1(\alpha\partial_2),
\qquad 1-\alpha\partial_i=\pi_i, \quad i=1,2.
\end{equation}
\end{Lemma}
\begin{Proof}
By the definitions
$$
(\pi_1\circ\pi_2)(u)=u-\alpha\partial_1\circ\partial_2(u).
$$
Then 
$$
\begin{array}{rl}
\alpha\partial_1\circ\partial_2(u)
&=u-\pi_1\circ\pi_2(u)=u-\pi_1(\pi_2(u))=u-\pi_1(u-\alpha\partial_2(u))
\\
&=u-(u-\alpha\partial_2(u))+\alpha\partial_1[u-\alpha\partial_2(u)]
=\alpha[\partial_1(u)+\partial_2(u)-\partial_1(\alpha\partial_2)(u)].
\end{array}
$$
Thus
$$
\partial_1\circ\partial_2=\partial_1+\partial_2-\partial_1(\alpha\partial_2).
$$
The result follows.
\end{Proof}
Let us give some known examples of divided difference operators.
\vspace{2mm}

(1) Let $R=K[[x]]$, and $\alpha(x), \phi(x)\in R$ where $\alpha(0)=0$.
 For $p\in R$ let
$$
\partial p(x)=\frac{p(x)-p(x-\alpha(x)\psi(x))}{\alpha(x)} . 
$$
Clearly $\partial$ always satisfies (\ref{eq12}), i.e. this is
a divided difference operator. We have 
$$
\pi p(x)=p(x)-\alpha\cdot\frac{p(x)-p(x-\alpha(x)\psi(x))}\alpha
=p(x-\alpha(x)\psi(x)).
$$
This is a usual translation of the argument.

(2) Let $R=K[[x,y]]$ and $\alpha=x-y$. Let
\begin{equation}\label{eq16}
\partial p(x,y)=\frac{p(x,y)-p(y,x)}{x-y}\,.
\end{equation}
Here $\pi(p(x,y))=p(y,x)$, i.e.  $\pi^2=1$ and  $\partial^2=0$.

For $R=K[[x_1,\dots,x_n]]$ one can define the operators $\partial_{ij}$ by
$$
\partial_{ij}p(x_1,\dots,x_n)
=\frac{p(\dots x_i\dots x_j\dots)-p(\dots x_j\dots x_i\dots)}{x_i-x_j}\,.
$$
These operators play an important part in the cohomology of the flag
manifolds. They have been used often (under the same name ``divided difference
operators''). Their algebra is very beautiful.

Let $f(x,y)\in K[[x,y]]$ be a series defining some formal group over
$R$.  Let $\alpha=f(x,\bar y)$ where $\bar y$ is the element inverse
to $y$ with respect to the formal group $f(x,y)$. Then the formula
\begin{equation}\label{eq17}
\partial p(x,y)=\frac{p(x,y)-p(y,x)}{f(x,\bar y)}
\end{equation}
defines a divided difference operator in $R=K[[x,y]]$. Here
$\pi(p(x,y))=p(y,x)$, $\pi^2=1$, and $\partial^2=\gamma\partial$
for some $\gamma\in K[[x,y]]$. For example, let $f(x,y)=x+y-axy$, $a\in K$.
Then $\partial^2=a\partial$. 

The operators given by (\ref{eq17}) with $f(x,y)$ being the formal group
of geometric cobordism elements have interesting applications in the
complex cobordism theory (see \cite{18}, \cite{19}).

The examples (1) and (2) are particular cases of the following general
construction. 

Let $R=K[[x_1,\dots,x_n]]$ and $\alpha(x)\in R$ with $\alpha(0)=0$.
Put
$$
\psi(x)=(\psi_1(x),\dots,\psi_n(x))
$$
where $x=(x_1,\dots,x_n)$ and 
$\psi_k(x)\in R$, $k=1,\dots,n$. 

Now put
$$
\partial p(x)=\frac1{\alpha(x)}
\bigl(p(x)-p(x-\alpha(x)\psi(x))\bigr)
$$
for any $p\in R$. We have
$$
\pi p(x)=p(x-\alpha(x)\psi(x)).
$$
For applications, the most interesting cases are when 
$\alpha(x)=\langle x,\xi\rangle=\sum^n_{i=1}x_i\xi^i$ and either

(i) $\psi(x)=\frac\xi{\langle\xi,\xi\rangle}$ where $\pi_\xi
p(x)=p\bigl(x-\langle
x,\xi\rangle\frac\xi{\langle\xi,\xi\rangle}\bigr)$ is a projector,
$\pi^2_\xi=\pi_\xi$, and $\partial=\partial_\xi$ is a division operator
by  $\langle x,\xi\rangle$, or

(ii) $\psi(x)=2\frac\xi{\langle\xi,\xi\rangle}$ where $\pi_\xi
p(x)=p\bigl(x-2\langle
x,\xi\rangle\frac\xi{\langle\xi,\xi\rangle}\bigr)$ satisfies 
the conditions $\pi^2_\xi=1$ and $\partial^2_\xi=0$.

These are the operators that are determined by reflections
with respect to the hyperplanes $V=\{x\mid\langle x,\xi\rangle=0\}$.
They were used to define the Dunkl's operators, see
\cite{20}.

Particular choices of a configuration space for the vectors
$\xi_1,\dots,\xi_k$ lead to interesting algebras of  the operators 
$\partial_{\xi_l}$, $l=1,\dots,k$.

(3) Let $\phi :R\to R$ be a ring homomorphism. According to \cite{21},
a {\it derivation of the algebra $R$} is an operator $\delta : R\to R$,
such that
\begin{equation}\label{eq18}
\delta(ab)=\phi(a)\delta(b)+\delta(a)b.
\end{equation}
\begin{Definition}\label{def4}{\rm 
An algebra over $R$ additively isomorphic to $R[t]$ is called 
an {\it Ore extension} if $R\subset R[t]$ and 
\begin{equation}\label{eq19}
ta=\phi(a)t+\delta(a).
\end{equation}
The operator $\phi$ from (\ref{eq19}) is called {\it 
Ore's $\phi$-derivation.}}
\end{Definition}
One easily proves the following lemma.
\begin{Lemma}\label{lemma10}
A divided difference operator $\partial$ satisfying the identity 
{\rm (\ref{eq12})}
$$
\partial(ab)=(\partial a)b+a(\partial b)-\alpha(\partial a)(\partial b)
$$
where $\alpha\in R$ is an Ore's $\phi$-derivation with
$\phi=\pi=1-\alpha\partial : R\to R$.
\end{Lemma}
Now we construct new $K$-linear product structures in the ring $R$.
\begin{Theorem}\label{theorem1}
Let $\alpha_1,\alpha_2\in R$ and let $\partial_1$, $\partial_2$ be
divided difference operators corresponding to the elements
$\alpha_1,\alpha_2$. Let $\pi_1 = 1-\alpha_1\partial_1$, $\pi_2 =
1-\alpha_2\partial_2$ be corresponding multiplicative projectors.
The operation
\begin{equation}\label{eq20}
\mu_1(x,y)=\pi_1(x)\pi_2(y)=xy-\alpha_1\partial_1(x)y
- -\alpha_2\partial_2(y)x+\alpha_1\alpha_2\partial_1(x)\partial_2(y)
\end{equation}
determines an associative product if and only if $\partial_1$ and
$\partial_2$ are division operators and the multiplicative projectors
$\pi_1$ and $\pi_2$ commute, i.e. $\pi_1\pi_2=\pi_2\pi_1$.  In
particular, each division operator $\partial$ determines the
associative product $\mu(x,y)=\pi(x)\pi(y)$ where
$\pi=1-\alpha\partial$. The product $\mu_1(x,y)$ is commutative if and
only if $\pi_1=\pi_2$.
\end{Theorem}
\begin{Proof}
Let 
$$
x\circ y=\mu_1(x,y)=\pi_1(x)\pi_2(y).
$$
Then 
\begin{equation}\label{eq21}
(x\circ y)\circ z=\pi_1(x\circ y)\pi_2(z)=\pi^2_1(x)\pi_1(\pi_2(y))\pi_2(z).
\end{equation}
Assume that the divided difference operators $\partial_1$ and
$\partial_2$ are division operators,
i.e. $\partial_1(\alpha_1)=\partial_2(\alpha_2)=1$. Thus
$\pi^2_1=\pi_1$, $\pi^2_2=\pi_2$, and the identity (\ref{eq21})
implies that the product $\mu_1$ is associative.

Conversely, let $\mu_1$ be associative. Then the condition
$$
(x\circ1)\circ1=x\circ(1\circ1)
$$
implies
$$
\pi^2_1(x)\pi_1(\pi_2(1))\pi_2(1)=\pi_1(x)\pi_2(\pi_1(1))\pi^2_2(1).
$$
We have  $\pi^2_1=\pi_1$ since $\pi(1)=1$ for any multiplicative operator 
$\pi$. Similarly, we have $\pi^2_2=\pi_2$. 

Now the condition $(1\circ y)\circ1=1\circ(y\circ1)$ implies that 
$\pi_1\pi_2=\pi_2\pi_1$. We use Corollary \ref{cor2} to conclude that
the operators $\partial_1$ and $\partial_2$ corresponding to the
projectors $\pi_1$ and $\pi_2$  are division operators.

Finally, if $\mu_1(x,y)$ is a commutative product then the condition
 $x\circ1=1\circ x$ implies that $\pi_1=\pi_2$. This proves Theorem
 \ref{theorem1}.
\end{Proof}
\begin{Theorem}\label{theorem2}
Let $\partial$ be a divided difference operator corresponding to an
element $\alpha\in R$ and $\beta\in R$. The
operation
\begin{equation}\label{eq22}
\mu_2(x,y)=xy+\beta\partial(x)\partial(y), \qquad \beta\in R,
\end{equation}
determines an associative product in $R$ if and only if one of the
following conditions is satisfied:
\begin{enumerate}
\item[{\rm (i)}] the operator $\partial$ is a division operator and 
$\pi(\beta)=0$  where $\pi=1-\alpha\partial$;
\item[{\rm (ii)}] the operator $\partial$ is not a division operator and
$\partial^2x\partial y=\partial x\partial^2y$ for any elements  $x,y\in R$.
\end{enumerate}
In particular, the condition {\rm (ii)} is satisfied if
$\partial^2=\gamma\partial$.

Furthermore, the condition {\rm (ii)} does not depend on a choice of
$\beta$, i.e. $\beta$ is a parameter of a deformation that connects the
original product structure in the ring $R$ to the  product
structure $\mu_2$.
\end{Theorem}
\begin{Proof}
Let
$$
x\circ y=\mu_1(x,y)=xy+\beta\partial x\partial y.
$$
Then 
$$
\begin{array}{rl}
(x\circ y)\circ z&=(xy+\beta\partial x\partial y)z
+\beta\partial(xy+\beta\partial x\partial y)\partial z
\\
&=xyz+\beta\bigl((\partial x\partial y)z+(\partial x)y\partial z
+x\partial y\partial z\bigr)
\\
&\qquad-\beta\alpha\partial x\partial y\partial z
+\beta\bigl(\partial\beta\partial x\partial y\partial z
+\beta\partial(\partial x\partial y)\partial z-\alpha\partial\beta
\partial(\partial x\partial y)\partial z\bigr).
\end{array}
$$
Thus the product $\mu_2$ is associative if and only if
\begin{equation}\label{eq23}
(\beta-\alpha\partial\beta)
\bigl(\partial(\partial x\partial y)\partial z
- -\partial x\partial(\partial y\partial z)\bigr)=0
\end{equation}
for all $x,y,z\in R$.

Let $\partial$ be a division operator, i.e. $\partial(\alpha)=1$.  Put
$x=\alpha^2$, $y=z=\alpha$. Then we have
$\partial(\alpha^2)=2\alpha\partial\alpha-\alpha(\partial\alpha)^2=\alpha$
whence $\partial^2(\alpha^2)=1$. 
Using $\partial(1)=0$
for the fixed values of
$x,y,z$ we see that
(\ref{eq23}) gives $\pi(\beta)=\beta-\alpha\partial\beta=0$. This
proves (i).

Now assume that $\partial$ is not a division operator. Then according to
Lemma \ref{lemma8}, $\pi(\beta)\ne0$ if $\beta\ne0$. The identity
$$
\partial(\partial x\partial y)=\partial^2x\partial y+\partial x\partial^2y
-\alpha\partial^2x\partial y
$$
shows that (\ref{eq23}) is equivalent to
\begin{equation}\label{eq24}
(\partial^2x\partial z-\partial z\partial^2z)
(\partial y-\alpha\partial^2y)=0
\end{equation}
for all $x,y,z\in R$. We notice that $\partial
y-\alpha\partial^2y=\pi\partial y$. The operator $\partial$ is not a
division operator by assumption thus $\pi\partial(y)=0$ is equivalent
to $\partial y=0$ for all $y\in R$. This is
impossible. It follows now that (\ref{eq24}) is equivalent to
$\partial^2x\partial z=\partial z\partial^2z$. This proves (ii) and
the theorem.
\end{Proof}
Now we give one more construction of an associative product which
comes from the complex cobordism theory (see \cite{22}).
\begin{Theorem}\label{theorem3}
Let $\Pi:R\to R$ and  $\delta:R\to R$ be linear operators such that
\begin{enumerate}
\item[{\rm (1)}] $\Pi^2=\Pi$, $\delta\Pi=\delta$,
\item[{\rm (2)}] $\delta(\Pi x\Pi y)=\delta x(\Pi y)
+(\Pi x)\delta y-\alpha(\delta x)(\delta y)$,
\item[{\rm (3)}] $\Pi((\Pi x)(\Pi y))=(\Pi x)(\Pi y)+\beta(\delta x)(\delta y)$
\end{enumerate}
for some elements $\alpha, \beta \in R$. Then the operation
$$
\mu(x,y)=\Pi((\Pi x)(\Pi y))=(\Pi x)(\Pi y)+\beta(\delta x)(\delta y)
$$
is an associative product structure in the ring $R$.
\end{Theorem}
\begin{Proof}
We denote $\mu(x,y)=x*y$. Then the definitions imply
$$
\begin{array}{rl}
(x*y)*z&= \ \Pi(x*y)\Pi z+\beta\delta(x*y)\delta z
\\
&= \ (\Pi(\Pi x\Pi y))\Pi z+\beta\delta(\Pi((\Pi x)(\Pi y)))\delta z
\\
&=\ \Pi x\Pi y\Pi z+\beta(\delta x)(\delta y)\Pi z
+\beta[\delta x(\Pi y)\delta z+(\Pi x)\delta y\delta z-\alpha(\delta x)
(\delta y)(\delta z)]
\\
&=\ \Pi x\Pi y\Pi z+\beta[\delta x\delta y\Pi z+\delta x\Pi y\delta z
+\Pi x\delta y\delta z]-\beta\alpha\delta x\delta y\delta z
\\
&= \  x*(y*z).
\end{array}
$$
The result follows.
\end{Proof}
\section{Examples from cobordism theory}
Here we use the above results to construct new product structures in
the complex cobordism theory.

First, we recall that all multiplicative operators in the algebra
$A^U=(\Lambda S)^\wedge$ are completely determined by their action on 
geometric cobordism elements $x\in U^2(X)$, i.e. they all are given as 
series 
\begin{equation}\label{eq25}
\phi(x)=x+\sum_{i\ge1}\phi_ix^{i+1}, \qquad \phi_i\in\Lambda.
\end{equation}
A multiplicative operator $\phi$ given by a series $\phi(x)$ acts on
an element $y\in U^*(X)$ as follows:
$$
\phi(y)=y+\sum_{\operatorname{deg}w>0}\phi_ws_w(y),
$$
where $\phi_w=\phi^{k_1}_1\dotsb\phi^{k_l}_l$ with 
$w=(k_1,\dots,k_l)$.
\begin{Lemma}\label{lemma11}
A series $D\phi(x)=x+\sum_{i\ge1}\phi_ix^{i+1}$ determines a
multiplicative projector $\phi\in A^U$ (i.e. $\phi^2=\phi$) if and
only if $\phi(\phi_i)=\phi_i+\sum\phi_ws_w(\phi_i)=0$ for all $i\ge1$.
\end{Lemma}
\begin{Proof}
Let $x\in U^2(X)$ be a geometric cobordism element. Then
$$
\phi(\phi(x))=\phi(x)+\sum\phi(\phi_i)\phi(x)^{i+1}.
$$
Thus $\phi(\phi(x))=\phi(x)$ if and only if
$\sum_{i\ge1}\phi(\phi_i)t^{i+1}=0$ in the ring $\Lambda[[t]]$ with
$t=\phi(x)$. The latter is equivalent to 
$\phi(\phi_i)=0$, $i\ge1$ and the result follows.
\end{Proof}
Let $x\in U^2(X)$ be a geometric cobordism element. A multiplicative
operator $\phi\in A^U$ with $\phi(x)=x+\sum_{i\ge1}\phi_ix^{i+1}$,
$\phi_i\in\Lambda$, is called {\it homogeneous} if $\deg \phi_i=-2i$
for all $i=1,2,\ldots$. The set of all homogeneous multiplicative
operators form a semigroup. Clearly multiplicative projectors 
belong to this semigroup. This semigroup is a fundamental object in the
complex cobordism theory. 

Let $m$ be an integer and $\Bbb Z_{(m)}=\Bbb Z[m^{-1}]$. The paper
\cite{16} contains the following result.
\begin{Lemma}\label{lemma12}
Let $\alpha\in\Lambda^{-2n}$ and $s_{(n)}\alpha=m\ne0$, $m\in\Bbb
Z$. An operator $\partial \in A^U$ is a division operator by an
element $\alpha$ in the localized complex cobordism theory
$U^*(X)\otimes\Bbb Z_{(m)}$ (so that $\pi=1-\alpha\partial$ is a
homogeneous multiplicative operator) if and only if the value of
$\partial$ on any geometric cobordism element $x$ is given by
$$
\partial x=\frac1mx^{n+1}+\sum_{i\ge1}a_ix^{n+i+1},
\qquad a_i\in\Lambda^{-2i}\otimes\Bbb Z_{(m)},
$$
where $(a_i)$ are free parameters of $\partial$.
\end{Lemma}
\begin{Proof}
According to the above algebraic results, it is enough to verify that
$\pi(\alpha)=0$. Indeed, we have: $\pi(x)=x-\alpha
x^n\bigl(\frac1mx+\sum_{i\ge1}a_ix^{i+1}\bigr)$ where $x$ is a
geometric cobordism element. Thus 
$$
\pi(y)=y-\frac1m \alpha\cdot s_{(n)}y+ 
(\mbox{terms containing $s_w(y)$ with $\deg w >2n$}).
$$
We conclude that $\pi(\alpha)=\alpha-\alpha\frac1ms_{(n)}(\alpha)=0$.
The result follows.
\end{Proof}
Lemma \ref{lemma12} may be used to describe all pairs
($\partial_1$,  $\partial_2$) of division operators such that
the corresponding multiplicative operators $\pi_1$ and
$\pi_2$ are homogeneous. Thus according to Theorem \ref{theorem1},
we obtain the set of multiplicative associative product structures
(in the complex cobordism theory)
of the form $\mu_1(x,y)=\pi_1(x)\pi_2(y)$.

Moreover, let $\partial$ be a division operator from Lemma
\ref{lemma12}, $\pi=1-\alpha\partial$, and $\beta \in \mbox{Ker}(\pi)$.
Then according to Theorem \ref{theorem2} (ii), we obtain an associative
product structure $\mu_2(x,y)=xy+\beta\partial x\partial y$.
 
The following result provides a construction of associative product
structures based on Theorem \ref{theorem2} (ii).
\begin{Lemma}\label{lemma13}
Let $\alpha\in\Lambda^{-2n}$ and $s_{(n)}(\alpha)=2n$.  Then any
multiplicative homogeneous operation $\pi\in A^U\otimes\Bbb Z_{(n)}$
given by
\begin{equation}\label{eq26}
\pi(x)=x-\alpha\partial x=\frac x{\root n\of{\vphantom{(}1+\alpha x^n}}
\end{equation}
on a geometric cobordism element $x\in U^2(X)$ determines the divided
difference operator $\partial$ such that $\partial^2=0$.
\end{Lemma}
\begin{Proof}
According to Lemma \ref{lemma6}, it is enough to verify that
$$
\pi(\alpha)=-\alpha \quad
\mbox{and}\quad \pi^2=1.
$$
We have
$$
\pi(x)=x-\frac1n\alpha x^{n+1}+O(x^{2n+1}).
$$
Thus 
$$
\pi(y)=y-\frac{1}{n} \alpha\cdot s_{(n)}y+ 
(\mbox{terms containing $s_w(y)$ with $\deg w >2n$}).
$$
Then it follows that
$\pi(\alpha)=\alpha-\frac1n\alpha\cdot2n=-\alpha$.  The action of
$\pi^2$ is determined by its value on a geometric cobordism element
since $\pi^2$ is a multiplicative operator. We have
$$
\pi^2(x)=\pi\biggl(\frac x{\root n\of{\vphantom{(}1+\alpha x^n}}\biggr)
=\frac{\pi(x)}{\root n\of{1+\pi(\alpha)\pi(x)^n}}
=\frac x{\root n\of{\vphantom{(}1+\alpha x^n}}\cdot
\frac1{\root n\of{1-\alpha\dfrac{x^n}{1+\alpha x^n}}}=x.
$$
Thus $\pi^2(x)=x$ whence $\pi^2=1$. The result follows.
\end{Proof}

\noindent
{\bf Example.}
Let $[\Bbb CP^m]$ be the cobordism class of $\Bbb CP^m$. Then
$s_{(m)}[\Bbb CP^m]=-(m+1)$. Hence we can choose $\alpha$ equal
to $-[\Bbb CP^{2n-1}]$ to obtain the divided difference operator
$\partial$ on $U^*(X)\otimes\Bbb Z_{(n)}$ (according to Lemma
\ref{lemma13}). Here $\partial$ satisfies $\partial^2=0$ which gives
a new associative product structure depending on free parameter
$\beta$. For instance, the complex projective plane $\Bbb CP^1$ gives
such a product on $U^*(X)$.
\vspace{2mm}

In conclusion, we describe particular example of an associative product
on $U^*(X)$ based on the construction from Theorem \ref{theorem3}.

Following Conner and Floyd, we define the additive projector
$$
\Pi: U^*(X)\to U^*(X)
$$
which is completely determined by the following action on the
complex cobordism classes.

Let $\xi\to\Bbb CP^1$ be a cononical line bundle and $M^{2n}$ 
a closed $U$-manifold. Let $\tau$ be its stable complex tangent
bundle. By definition, $\Pi([M^{2n}])$ is a cobordism class of the
manifold $i : \widehat M^{2n}\subset M^{2n}\times\Bbb CP^1$ with the
normal bundle $i^*((\det\tau)\otimes\xi)$ where $\det\tau$ is the
determinant bundle.

The action of the projector $\Pi\in A^U$ on the Thom class $x_n$ (of a
 complex $n$-dimensional bundle $\eta\to X$) is given by
$$
\Pi(x_n)=x_n+\sum_{i\ge2}\alpha_{i1}\partial_ix_n
$$
where $\alpha_{i1}\in\Lambda^{-2i}$ is the coefficient with  $x^iy$ 
in the formal group $f(x,y)$ of geometric cobordism elements $x,y$.
Here  $\partial_i\in A^U$ are the projectors acting on the Thom class
$x_n$ as
$$
\partial_ix_n=x_nc_1(\overline{\det\eta})
$$
where $c_1$ is the first Chern class in complex cobordism and
``$\overline{\phantom{n}}$'' is the complex conjugation.

Let $\delta=\partial_1\in A^U$. Then one can verify that the pair
($\Pi,\delta$) satisfies the conditions of Theorem \ref{theorem3} (see
\cite{22}, \cite{23} for details). Thus according to Theorem
\ref{theorem3}, we obtain the associative product given by
$$
x*y=\Pi(\Pi x\Pi y).
$$
This product structure is crucial to describe a ring structure of the
cobordism ring of $SU$-manifolds (these manifolds are called sometimes
Calabi-Yau manifolds).

\end{document}